\documentclass[12pt]{article}
\usepackage{theorem}
\newtheorem{lemma}{Lemma}

\newtheorem{theorem}[lemma]{Theorem}
\newtheorem{corollary}[lemma]{Corollary}
{\theorembodyfont{\upshape}}
{\theorembodyfont{\upshape}}
{\theorembodyfont{\upshape}}
{\theorembodyfont{\upshape}}

\newcommand{\R}{{\bf R}}

\newcommand{\rme}{{\rm e}}
\newcommand{\rmd}{{\rm d}}

\newcommand{\alp}{\alpha}
\newcommand{\bet}{\beta}
\newcommand{\gam}{\gamma}
\newcommand{\lam}{\lambda}

\newcommand{\eps}{\varepsilon}

\newcommand{\Ome}{\Omega}
\newcommand{\lap}{{\Delta}}
\newcommand{\gra}{\nabla}

\newcommand{\tr}{{\rm tr}}
\newcommand{\Dom}{{\rm Dom}}

\newcommand{\norm}{\Vert}

\newcommand{\Proof}{\underbar{Proof}{\hskip 0.1in}}

\newcommand{\dist}{{\rm dist}}
\newcommand{\Schrodinger}{Schr\"odinger }
\newcommand{\capac}{{\rm cap}}

\newcommand{\pr}{\prime}
\newcommand{\nolabel}{\nonumber}
\title{ A REVIEW OF HARDY INEQUALITIES}
\author{E.B. Davies}
\date{August 1998}
\begin{document}
\maketitle
\begin{abstract}
We review the literature concerning the Hardy inequality 
for regions in Euclidean space and in manifolds, 
concentrating on the best constants. We also give 
applications of these inequalities to boundary decay and spectral 
approximation. 
\par
keywords: Hardy inequality, boundary decay, Laplacian, 
elliptic operators, spectral theory, eigenfunctions, 
spectral convergence.
\end{abstract}
\section{Introduction}
\par
Let $H$ be a non-negative second order elliptic operator 
acting in $L^{2}(U)$ subject to Dirichlet boundary 
conditions, where $U$ is a region in $\R^{N}$ or 
in a Riemannian manifold. Also let $d$ be a positive 
function on $U$ which is continuous and satisfies $|\gra 
d |\leq 1$. Traditionally one takes $d(x)$ to be the distance 
of $x\in U$ from the boundary $\partial 
U$, but another possibility is that $d(x)$ is the distance 
from any closed subset of $M\backslash U$ if $U$ is embedded 
in some larger Riemannian manifold $M$. 

We say that $H$ satisfies a weak Hardy inequality with 
respect to $d$ if there exists a
constant $c>0$ and a constant 
$a\geq 0$ such that
\begin{equation}
\int_{U}\frac{|f|^{2}}{d^{2}}\leq c^{2}\left( Q(f) +a\norm f 
\norm^{2 }\right)	
	\label{1}
\end{equation}
is valid for all $f\in C_{c}^{\infty}(U)$, and hence for all $f$ in 
the domain of the quadratic 
form $Q$ of $H$. The infimum of all possible $c$ in 
(\ref{1}) is then 
called the weak Hardy constant. We say that $H$ satisfies a 
strong Hardy inequality if (\ref{1}) holds with $a=0$, in 
which case the minimum possible $c$ is called the strong 
Hardy constant.

There are also $L^{p}$ and higher order analogues of the above notion, which 
we mention briefly later in this review.

In section 2 we describe the method of geodesic integrals 
for proving Hardy inequalities in higher dimensions. 
Section 3 describes a method ultimately due to Jacobi, 
while Section 4 gives various miscellaneous results.

We then turn to the applications of the HI to the 
proof of boundary decay. It was shown 
in \cite{D1} that Hardy's inequality can be used to prove 
the $L^{2}$ boundary decay of eigenfunctions without any 
further assumptions. This in turn leads to the possibility 
of controlling the rate of convergence of the eigenvalues 
when the region $U$ is approximated by a family $U_{\eps}$ 
of slightly smaller regions. Very recently progress has 
been made on this problem, \cite{D4}, and we are able to announce 
bounds on the rate of convergence which are sharp 
in a certain sense.

Our main results on boundary decay, Theorems 11 and 12, may 
be regarded as $L^{2}$ analogues of much stronger pointwise bounds on 
eigenfunctions given in \cite{Ba2,CZ,LP}. Note however that our 
bounds depend only on the validity of (\ref{1}), hold for 
all functions in the domains of 
the operators, not just for eigenfunctions, and have rather 
precise constants. 

If we abandon interest in the 
precise value of the constant, and choose $d$ to be the 
Euclidean distance from an arbitrary point of $U$, then it 
may be seen that our results are related to Morrey space 
estimates. These have been of considerable importance in the 
theory of elliptic operators, and recently in the proof of 
heat kernel bounds, and we refer the reader to 
\cite{Au,AT,Gi} for further details. 
\section{Geodesic integrals}

The first method which we describe depends upon 
the one-dimensional case, which is the only one Hardy actually 
studied. We refer to \cite{OK} for an exhaustive study, which 
involves generalizations to the variable 
coefficient case of the original formula
\[
\int_{0}^{\infty}\frac{|f(x)|^{2}}{x^{2}} \rmd x \leq 
4\int_{0}^{\infty}|f^{\pr}(x)|^{2}\rmd x
\]
valid for all $f\in C_{c}^{\infty}(0,\infty)$, and hence for 
all $f\in W^{1,2}_{0}(0,\infty)$.

Let $H:= -\lap_{DIR}$
in the Hilbert space $L^{2}(U)$ where $U$ is a bounded 
region in $\R^{N}$. For every unit vector $u\in S^{N-1}$ 
and $x\in U$ let
\[
 d_{u}(x):=\min \{ |t| :x+tu\notin U\}
\]
if the set of such $t$ is non-empty, and put $d_{u}(x):=+\infty$ 
otherwise.
We define the (harmonic) mean distance of $x$ from $\partial 
U$ by
\begin{equation}
m(x)^{-2}:=|S^{N-1}|^{-1}\int_{S^{N-1}}
 d_{u}(x)^{-2} \rmd S(u).                \label{mdist}  
\end{equation}
It is easy to prove that $d(x)\leq m(x)$ for all $x\in U$.
\begin{lemma}
We have
\[
\frac{N}{4 m^{2}} \leq H
\]
in the sense of quadratic forms. If $\lam_{1}$ is the 
smallest eigenvalue of $H$ then 
\[
\lam_{1}\geq \frac{N}{4\mu^{2}}
\]
where the quasi-inradius $\mu$ of $U$ is defined by
\[
\mu:=\sup \{ m(x):x\in U\} .
\]
\end{lemma}

\Proof See \cite{D3} or \cite[Th. 1.5.3]{HKST}.

Applications of the above lemma depend on making 
assumptions on $U$ which enable one to bound $m(x)$ above 
by some multiple of $d(x)$. The first of these is folklore 
and seems not to 
have been written down explicitly until very recently; see 
\cite{MMP,MS} and the next section for alternative proofs.

\begin{theorem}
If $U$ is a convex subset of $\R^{N}$ then 
\[
\frac{1}{d^{2}} \leq 4 H
\]
in the sense of quadratic forms.
\end{theorem}

\Proof If $a$ is the point of $\partial U$ closest to 
$x$ then we can obtain the relevant upper bound of $m(x)$ by 
computing an appropriate integral over the supporting hyperplane 
at $a$. See \cite[Exercise 5.7]{STDO}.

The following lemma is typical of a variety of methods of 
obtaining crude upper bounds on $m(x)$. The hypothesis is 
valid not only for regions with Lipschitz boundaries, but 
also for a variety of regions with fractal boundaries, 
such as the Koch snowflake region in $\R^{2}$. 
\begin{lemma}
Suppose that there is a constant $k$ such that 
for each $a\in \partial U$ and each $\alp >0$ 
there exists a ball $B$ disjoint from $U$ with centre $b$ 
and radius $\bet\geq k\alp$, where $|b-a|=\alp$. Then there 
exists constants $c_{0}, c_{1}$ such that $m(x)\leq c_{0}d(x)$ and 
hence
\[
\frac{1}{ d^{2}} \leq c_{1}H
\]
in the sense of quadratic forms.
\end{lemma}

\Proof See \cite{D3}, \cite[Th. 1.5.4]{HKST} and \cite[Th. 
3]{A}.

The condition of Lemma 3 is not satisfied for regions 
satisfying a uniform exterior power-like cusp condition. In 
such cases one may prove a modified Hardy inequality using 
Lemma 1, namely
\begin{equation}
\int_{U}\frac{|f|^{2}}{d^{\gam}}\leq c^{2}\left( Q(f) +a\norm f 
\norm^{2 }\right)	\nolabel
\end{equation}
for some $0<\gam <2$; see \cite[p 369]{DS}. See also 
\cite[Th. 3.2, 3.3]{DL} where a similar situation arises 
for locally Euclidean manifolds with fractal boundaries.

A procedure closely related to the idea of this section 
was developed for regions in 
Riemannian manifolds independently by Croke and Derdzinski, \cite{CD}, and 
Donnelly, \cite{Don}. The integrals over straight 
lines were replaced by integrals over geodesics, so the 
formulation involves the geodesic flow on the unit 
sphere bundle of the manifold. However, both papers are 
concerned with obtaining lower bounds on the 
bottom eigenvalue, much as in Lemma 1, rather than Hardy's 
inequality.

We mention in passing that there is no requirement that one 
should assign equal weights to every direction in Euclidean 
space. In some cases one obtains a better constant in the 
Hardy inequality by taking 
an average over a few directions which are well adapted to 
the region in question. 

\section{The Classical Method}

The following method goes back to Jacobi, 
and was used by Barta and Kasue to obtain lower bounds on 
the first eigenvalue, \cite{Ba,Ka}. It is the easy half of a
theorem of Allegretto, Moss and Piepenbrink  
characterising the bottom of the spectrum of a \Schrodinger 
operator in terms of the existence of positive 
distributional solutions of 
the eigenvalue equation, \cite[p. 23]{CFKS}. Assume that 
\[
Hf(x):=-\sum\frac{\partial}{\partial x_{i}}\left\{ 
a_{i,j}(x)\frac{\partial f}{\partial x_{j}} \right\}
\]
where $a(x)$ is a non-negative $C^{1}$ real symmetric matrix-valued 
function and $f\in C_{c}^{2}(U)$. Then $H$ is a 
non-negative symmetric operator and we can use the same symbol 
to denote its Friedrichs extension.

\begin{lemma}
Let $\phi$ be a positive $C^{2}$ function on $U$ and let $V$ 
be a continuous function on $U$ such that
\[
-\sum\frac{\partial}{\partial x_{i}}\left\{ 
a_{i,j}(x)\frac{\partial \phi}{\partial x_{j}} \right\} 
\geq V\phi .
\]
Then we have
\[
H\geq V
\]
in the sense of quadratic forms.
\end{lemma}

\Proof See \cite[Th. 4.2.1]{HKST}. 

The conditions of the above lemma can be weakened to allow a 
distributional inequality.

\underbar{Second proof of Theorem 2} If we put 
$\phi:=d^{1/2}$ and use the fact that $\lap d \leq 0$ for 
any convex set $U$ then the result follows immediately from 
the last lemma. 

The method of this section can be extended to 
Riemannian manifolds without difficulty.  We refer to  \cite{Ca,DH,Ow} 
for a variety of Hardy and Rellich type inequalities with 
explicit constants in Riemannian manifolds obtained in this 
manner. The following theorem is only one of a range of related 
results due to Brezis and Marcus, \cite{BM}. In particular they find explicit 
bounds on the minimum possible negative value of $a$ in the theorem 
when $U$ is convex.

\begin{theorem} 
If $\, U\subseteq \R^{N}$ is bounded 
with a $C^{2}$ boundary and $H:=-\lap_{DIR}$ in 
$L^{2}(U)$ then there exists $a\in\R$ such that
\begin{equation}
d^{-2}\leq 4(H+a)     \label{bm}
\end{equation}
in the sense of quadratic forms. If $U$ is convex 
then (\ref{bm}) holds for certain $a<0$.
\end{theorem}

\Proof  Let $\phi$ be a positive $C^{2}$ function on $U$ 
such that $\phi(x)=d(x)^{1/2} -d(x)$ 
for all $x$ close enough to $\partial U$. The first 
statement of the theorem 
follows by applying Lemma 4 to $\phi$.

There are various other improvements of the strong Hardy 
inequality of which we mention just two. For a definitive 
treatment of the one-dimensional theory see \cite{OK}.

\begin{theorem}
If $U:=\{ x\in\R^{N}:x_{N}>0\}$ where $N>1$ then
\[
\int_{U}\left\{  
\frac{1}{x_{N}^{2}} 
+\frac{1}{4x_{N}(x_{N}^{2}+x_{N-1}^{2})^{1/2}}
\right\} |f|^{2} \rmd x\leq 
4\int_{U}|\gra f|^{2}\rmd x
\]
for all $f\in C_{c}^{\infty}(U)$.
\end{theorem}

\Proof See \cite[Sect. 2.1.6]{Maz}.

\begin{theorem}
If $U:=(0,a)$ then
\[
\int_{U} \frac{a^{2}|f|^{2}}{x^{2}(a-x)^{2}}\rmd x\leq 
4\int_{U}| f^{\pr}|^{2}\rmd x
\]
for all $f\in C_{c}^{\infty}(U)$.
\end{theorem}

\Proof Put $\phi(x):=x^{1/2}(a-x)^{1/2}$ in Lemma 4.

\section{Capacity-based methods}

In this section we mention a few of the very general theorems 
which involve the use of capacity arguments. These have 
been developed in an $L^{p}$ context, but we only treat the 
case $p=2$. If $K$ is a compact subset of $U\subseteq 
\R^{N}$ we define its relative capacity by
\[
\capac (K,U):=\inf\left\{
\int_{U}|\gra f|^{2}: f\in C_{c}^{\infty}(U) \mbox{ {\rm 
and} } 
f\vert_{K}\geq 1
  \right\}.
\]
It is particularly appropriate in this conference 
to mention one version of the most quantitatively precise 
theorems of this type, due to Professor Maz'ya.

\begin{theorem}
If $\mu$ is a positive measure on $U$ and 
\[
\mu(K)\leq \bet \,\capac(K,U)
\]
for all compact subsets $K$ of $U$, then
\[
\int_{U}|f|^{2}\rmd \mu \leq 4\bet \int_{U}|\gra f|^{2}
\]
for all $f\in C_{c}^{\infty}(U)$. Conversely the second 
inequality implies 
\[
\mu(K)\leq 4\bet \,\capac(K,U)
\]
for all compact subsets $K$ of $U$.
\end{theorem}

\Proof See \cite[p.113]{Maz}.

Our next results are taken from a paper of Ancona, 
\cite{A}. We say that $U$ is uniformly $\lap$-regular if 
for all $x\in\partial U$ and all $r>0$
the harmonic measure $w$ of $U\cap \partial B(x,r)$ in 
$U\cap B(x,r)$ satisfies
$w\leq 1-\bet$ on $U\cap \partial B(x,r/2)$, for some 
constant $\bet\in(0,1)$ independent of $x,r$. If $N\geq 3$ 
this is equivalent to the uniform capacitary density 
condition that there 
exists a constant $\alp>0$ such that
\[
\capac(B(x,r)\backslash U) \geq \alp r^{N-2}
\]
for all $x\in\partial U$ and all $r>0$.

\begin{theorem} \label{CapTh}
If $N\geq 2$ and $U\subseteq \R^{N}$ 
is uniformly $\lap$-regular then $U$ 
satisfies a strong Hardy inequality with respect to the 
Laplace operator. If $N=2$ then the converse is also true.
\end{theorem}

Although \cite{A} does not provide sharp information 
about the size of the strong Hardy constant, it contains many 
more results than we have indicated above. An $L^{p}$ 
converse of Theorem \ref{CapTh} for $N=p>2$
may be found in \cite{Lew}, 
using an appropriate $L^{p}$ Riesz capacity.

\section{Miscellaneous results}

The weak Hardy constant $c$ as defined in Section 1 
was proved in \cite{D2} to be local in the sense that it is 
the maximum value of a certain upper semi-continuous function on 
the boundary, whose value at each point depends only on the 
geometry of the boundary around that point. Various methods 
of evaluating this function at different types of boundary 
point are described in \cite{D2}. 

For the remainder of this section we 
assume that $H:=-\lap_{DIR}$. The strong Hardy constant is a global 
invariant of $U$.  It equals $2$ for any convex set, 
but the condition of convexity is not necessary for 
this conclusion. Let 
\[
U_{\bet}:=\{ r\rme^{i\theta}:0<r<1 \,\,{\rm and}\,\, 0<\theta <\beta\}.
\]
Then $U_{\bet}$ has strong Hardy constant $2$ if and only if 
the internal angle $\bet$ is less than  or equal to a certain critical 
value $\bet_{c}\sim 4.856$ radians, \cite{D2}. For larger 
$\bet$ the strong and weak Hardy constants are larger than 
$2$. Similar conclusions hold for other plane regions with 
piecewise smooth boundaries.

If $U$ is a simply connected region in 
$\R^{2}$ then $U$ has strong Hardy constant at most $4$ by \cite{A}, 
\cite[Th. 1.5.10]{HKST}. The proof of this result depends 
upon a fact from analytic function theory, namely Koebe's 
one-quarter theorem.

There is an interesting connection between the possible 
constants in the strong Hardy inequality and the Minkowski 
dimension of the boundary, \cite{DM}. In two dimensions there is also a 
relationship with hyperbolic geometry, which we do not pursue.
We say that the boundary $\partial U$ has interior Minkowski dimension  
$\alp$ if there exist positive constants $k_{1}$ and $k_{2}$ 
such that
\[
k_{1}\eps^{N-\alp}\leq 
|\{ x\in U:\dist (x,\partial U)<\eps\}| \leq
 k_{2}\eps^{N-\alp}
\]
for all $\eps >0$. The following theorem is adapted from 
\cite[Th. 3.3]{DM}. We allow $\alp <N-1$ because the theorem 
is applicable in manifolds, for example if $U$ is obtained 
by removing a compact set $K$ from a sphere endowed with the 
standard metric.

\begin{theorem}
If $\partial U$ has interior Minkowski dimension  
$\alp>N-2$, then the strong Hardy 
constant of $U$ with respect to the Laplacian satisfies
\[
c(2+\alp -N)\geq 2.
\]
\end{theorem}

In most of the above lemmas we have restricted attention to 
Hardy inequalities in $L^{2}$. In fact many of the results 
have been extended to $L^{p}$ with sharp constants; see 
\cite{MMP,MS} for the proofs of the following two theorems.

\begin{theorem}
Let\[
c^{-p}:=\inf\left\{ \frac{ \int_{U}|\gra f|^{p}}{\int_{U}|f/d|^{p}}
:f\in W^{1,p}_{0}(U)
\right\}
\]
where $1<p<\infty$. If $\partial U$ is smooth then
$c\geq p/(p-1)$. If in addition $p=2$ then $c>2$ if and only 
if the infimum is achieved by some $f\in W^{1,2}_{0}(U)$.
\end{theorem}

\begin{theorem}
If $U$ is a convex set in $\R^{N}$ and $1<p<\infty$ then
\[
\int_{U} \frac{|f|^{p}}{d^{p}} \leq \left( 
\frac{p}{p-1}\right)^{p}\int_{U}|\gra f |^{p}
\]
for all $f\in W_{0}^{1,p}(U)$.
\end{theorem}

We refer to  \cite{A,BM,Ca,Lew,MMP,Maz,OK,W1,W2} for further $L^{p}$ 
results, since they do not yet have such direct 
consequences for spectral theory. We refer to 
\cite{DH, Ow} for the analogues for 
higher order operators, known as Rellich inequalities, and 
to \cite{Cia} for analogues in Orlicz spaces.

We describe some trace inequalities in 
\cite{D5} which may be proved using Hardy's 
inequality. We assume that $H:=-\lap_{DIR}$ acting in 
$L^{2}(U)$ where $U$ is a region in $\R^{N}$. The theorems
are only of interest when $U$ has infinite volume.

\begin{theorem}
We have
\[
\tr [\rme^{-Ht}] \leq (2\pi t)^{-N/2}
\int_{U}\rme^{-Nt/8m(x)^{2}}\rmd^{N}x
\]
for all $t>0$, where $m$ is defined by (\ref{mdist}).
\end{theorem}

\begin{theorem}
If $U$ satisfies the regularity condition
\[
d(x)\leq m(x)\leq bd(x)
\]
for all $x\in U$ then
\[
2^{-N}(2\pi t)^{-N/2}
\int_{U}\rme^{-8\pi^{2}N^{2}t/d(x)^{2}}\rmd^{N}x
\leq\tr [\rme^{-Ht}] 
\leq (2\pi t)^{-N/2}
\int_{U}\rme^{-Nt/8b^{2}d(x)^{2}}\rmd^{N}x
\]
for all $t>0$. Hence
\[
\tr [\rme^{-Ht}]<\infty
\]
for all $t>0$ if and only if
\[
\int_{U}\rme^{-t/d(x)^{2}}\rmd^{N}x <\infty
\]
for all $t>0$.
\end{theorem}

\section{Boundary estimates}

The size of the constant $s$ in an inequality of the form 
\begin{equation}
\int_{U}\frac{|f|^{2}}{d^{s}} <\infty   \label{power}
\end{equation}
conveys information about the behaviour of the function $f$ 
near the boundary of $U$. 
We conjecture that it is not possible to have $s>2$ in the inequality 
for any region $U$ 
if we only assume that $f\in\Dom(Q)$, where $Q$ is the 
quadratic form associated with a uniformly elliptic second 
order operator $H$ acting in $L^{2}(U)$ subject to 
Dirichlet boundary conditions. However, if 
we make stronger assumptions on $f$ then one may be able to 
prove (\ref{power}) for a larger value of 
$s$. The first paper with results of 
this type was \cite{EHK}, where it was assumed that $f$ was 
an eigenfunction of $H$. Subsequently \cite{D1} obtained 
better bounds for all $f\in\Dom(H)$, assuming only the Hardy 
inequality. 

Although we have concentrated on $L^{2}$ boundary 
estimates, there is a substantial literature on pointwise 
decay of eigenfunctions and their gradients at the boundary.
Bounds of the type
\begin{equation}
|\phi_{n}(x)|  \leq   c_{n} \phi_{1}(x)  \label{phi1}
\end{equation} 
are immediate consequences of intrinsic ultracontractivity 
(IU), 
\cite{DS,HKST}, in which a major ingredient of the proof 
is the existence of an inequality
\begin{equation}
\phi_{1}(x) \geq ad(x)^{\alp}           \label{lower}
\end{equation}
for some positive constants $a$ and $\alp$. The proof of 
(\ref{lower}) depends in turn upon the Harnack inequality 
and a boundary accessibility 
property. The BAP was proved in \cite{DS,HKST} for 
Lipschitz domains, but Ancona and Simon commented 
that it holds under a suitable 
twisted interior cone condition, i.e. for John domains, 
\cite[p 98]{Green}. Finally Banuelos gave a 
detailed analysis of the relationship between (\ref{phi1}), 
IU, John domains, Holder domains, NTA domains. etc. in \cite{Ba1}.
 
Pointwise bounds on the gradients of the
eigenfunctions $\phi_{n}$ of $-\lap_{DIR}$ and of  
\Schrodinger operators with potentials in restricted Kato 
classes acting in 
$L^{2}(U)$ are proved in 
\cite{Caf,CZ,LP,Ba2,BP} in steadily increasing generality. The 
best upper bound is for IU domains and is in \cite{Ba2}, 
while the best lower bound is for Lipschitz domains and is 
in \cite{BP}. The inequalities are of the form
\begin{eqnarray*}	
	|\gra\phi_{n}(x)| & \leq  & c_{n} \phi_{1}(x) /d(x)\\
	|\gra\phi_{1}(x)| & \geq  & c_{1} \phi_{1}(x) /d(x),
\end{eqnarray*}
the latter being for $x$ close enough to the boundary.

One may also obtain upper bounds of the form
\[
|\phi_{n}(x)|\leq c_{n}d(x)^{\bet}
\]
for explicit but non-optimal constants $c_{n},\bet$ which 
depend only on the eigenvalue $\lam_{n}$, the dimension and 
the constant $\alp$ in the uniform capacitary density 
inequality, \cite{BB}. For an open simply connected region 
in $\R^{2}$ the bound
\[
|\phi_{n}(x)|\leq c_{n}d(x)^{1/2}
\]
is proved in \cite{Ba2,BB,P2}; in this case the power $1/2$ is 
sharp.

We finally present some new results on $L^{2}$ boundary 
decay, taken from \cite{D4}. Let $U$ be 
a bounded region in $\R^{N}$ and let 
$H:=-\lap_{DIR}$ acting in $L^{2}(U)$. Let $d(x)$ denote the 
distance of $x$ from some closed subset of $\R^{N}\backslash 
U$. We make no assumptions on the boundary $\partial U$ 
apart from the validity of (\ref{1}) for certain values 
of $c\geq 2$ and $a\geq 0$. We are then able to draw the 
following conclusions 
about the boundary decay of functions in the domain of $H$. 
We have proved in \cite{D4} that the powers of $\eps$ in these theorems 
are sharp, and conjecture that the constant $c_{0}$ is also 
sharp. 
Analogues of the theorems for uniformly elliptic operators 
in divergence form are proved in \cite{D4}.

\begin{theorem}\label{5A}
If $f\in\Dom(H)$ and $\eps >0$ then
\[
\int_{\{x:d(x) <\eps  \}}|f|^{2}\leq 
c_{0}\eps^{2+2/c}\norm (H+a)f\norm_{2} 
\norm(H+a)^{1/c}f\norm_{2}
\]
where
\[
c_{0}:=c^{2+2/c}.
\]
\end{theorem}

\begin{theorem} \label{5B}
If $f\in\Dom (H)$ and $\eps >0$ then
\[
\int_{\{x:d(x) <\eps\}}|\gra f|^{2}
\leq c_{1}\eps^{2/c}\norm 
(H+a)f\norm_{2}\norm(H+a)^{1/c}f\norm_{2}.
\]
where
\[
c_{1}:=c^{2/c}+c^{2/c}(1+c)^{2+2/c}.
\]
\end{theorem}

\begin{corollary}
If $Hf=\lam f$ for some $\lam > 0$, $\norm f \norm_{2}=1$ and $\eps >0$ then
\[
\int_{\{x:d(x) <\eps  \}}|f|^{2}\leq 
c_{0}\eps^{2+2/c}(\lam +a)^{1+1/c}
\]
and
\[
\int_{\{x:d(x) <\eps\}}|\gra f|^{2}
\leq c_{1}\eps^{2/c}(\lam +a)^{1+1/c}.
\]
\end{corollary}

\begin{corollary}
If $U$ is a 
simply connected proper subregion of $\R^{2}$ then 
\[
\int_{\{x:d(x) <\eps  \}}|f|^{2}\leq 
32\eps^{5/2}\norm Hf\norm_{2} 
\norm H^{1/4}f\norm_{2}
\]
and
\[
\int_{\{x:d(x) <\eps  \}}|\gra f|^{2}\leq 
114\eps^{1/2}\norm Hf\norm_{2} 
\norm H^{1/4}f\norm_{2}
\]
for all $f\in\Dom(H)$ and $\eps >0$.
\end{corollary}
\Proof We put $c=4$ and $a=0$ in Theorems \ref{5A} and 
\ref{5B}.

We use the results above to consider the 
effect on the spectrum of $H:=-\lap_{DIR}$ of replacing the 
bounded region $U$ 
by a slightly smaller region $U_{\eps}$ such that
\[
\{ x\in U:d(x) >\eps\} \subseteq U_{\eps}\subseteq U.
\]
If $\lam_{n}(U_{\eps})$ denote the eigenvalues of the 
operator $H_{\eps}$ defined by restricting $H$ to $L^{2}(U_{\eps})$ 
where we again impose Dirichlet boundary conditions, then 
variational arguments imply that $\lam_{n}(U)\leq 
\lam_{n}(U_{\eps})$
for all $n$ and $\eps >0$. Our theorem below provides
quantitative estimates of the difference, again only 
assuming (\ref{1}). The first version in \cite{D1} did not 
obtain what we believe to be the sharp power of $\eps$ given 
below. Pang, \cite{P2}, obtained the result of Theorem 
\ref{5E} for 
$n=1$ for simply connected
plane regions by a method involving conformal mappings, 
improving his own earlier 
results in \cite{P1}. See \cite{D4} for a more general version 
of the theorem below, and its proof.

\begin{theorem}\label{5E}
There exist constants $c_{n}$ for all positive integers 
$n$ such that 
\[
\lam_{n}(U)\leq \lam_{n}(U_{\eps})\leq 
\lam_{n}(U)+c_{n}\eps^{2/c}.
\]
\end{theorem}

We finally mention that there is extensive literature which 
compares $\lam_{n}(U)$ with $\lam_{n}(U\backslash K)$, 
where $K$ is a compact subset of $U$ which has a small capacity 
in a suitable sense; we believe that \cite{RT} is one of the earliest 
contributions to this subject, often known as the crushed 
ice problem. See \cite{DMN} for a survey, including an 
explicit asymptotic formula for the difference of the 
eigenvalues in the limit 
of small ${\rm Cap}(K)$ and also estimates of the 
difference for $n=1$, both proved in the abstract setting of 
regular Dirichlet forms. See also \cite{No}, where 
estimates of the difference for $n=1$ in terms of an appropriate 
definition of capacity are obtained in an abstract contect 
applicable to higher order elliptic operators.
\par 
\vskip 0.3in
{\bf Acknowledgments } I should like to thank R Banuelos, 
H Brezis, V G Maz'ya and M Pang for comments on an 
early draft of this paper.
\par
\vfil
\eject

\par
\vskip 0.2in
Department of Mathematics \newline
King's College \newline
Strand \newline
London WC2R 2LS \newline
England\newline  
e-mail: E.Brian.Davies@kcl.ac.uk
\vfil
\end{document}